\documentclass[12pt,a4paper,oneside]{amsart}
\usepackage[a4paper]{geometry}
\usepackage{mathptmx} 
\DeclareMathAlphabet{\mathcal}{OMS}{cmsy}{m}{n} 

\usepackage{amssymb,perpage,stmaryrd}

\newtheorem{theorem}{Theorem}[section]
\newtheorem{proposition}[theorem]{Proposition}
\newtheorem{lemma}[theorem]{Lemma}
\newtheorem{corollary}[theorem]{Corollary}
\newtheorem{conjecture}[theorem]{Conjecture}
\newtheorem*{definition}{Definition}

\newcommand{\set}[2]{\ensuremath{\{ #1 \>|\> #2 \}}}

\def\liebrack{\ensuremath{[\,\cdot\, , \cdot\,]}}

\MakePerPage[2]{footnote}

\begin{document}

\title{Commutative $2$-cocycles on Lie algebras}

\author{Askar Dzhumadil'daev}
\address{Kazakh-British Technical University, Tole bi 59, Almaty 050000, Kazakhstan}
\email{dzhuma@hotmail.com}

\author{Pasha Zusmanovich}
\address{Reykjav\'ik Academy, Iceland}
\email{pasha.zusmanovich@gmail.com}

\date{last minor revision September 9, 2017}
\thanks{J. Algebra \textbf{324} (2010), N4, 732--748; \textsf{arXiv:0907.4780}}

\begin{abstract}
On Lie algebras, we study commutative $2$-cocycles, i.e., 
symmetric bilinear forms satisfying the usual cocycle equation.
We note their relationship with antiderivations and compute them for some classes of 
Lie algebras, including finite-dimensional semisimple, current and Kac-Moody algebras.
\end{abstract}

\maketitle

\section*{Introduction}

A \textit{commutative $2$-cocycle} on a Lie algebra $L$ over a field $K$ is a symmetric bilinear
form $\varphi: L \times L \to K$ satisfying the cocycle equation 
\begin{equation}\label{cocycle}
\varphi([x,y],z) + \varphi([z,x],y) + \varphi([y,z],x) = 0
\end{equation}
for any $x,y,z\in L$.
Commutative $2$-cocycles appear at least in two different contexts.
First, they appear in the study of 
nonassociative algebras satisfying certain skew-symmetric identities.
It turns out that all skew-symmetric identities of degree $3$ reduce to a number
of identities, among which
$$
[a,b]c + [b,c]a + [c,a]b = 0
$$
and
$$
a[b,c] + b[c,a] + c[a,b] = 0
$$
play a prominent role (here multiplication in an algebra is denoted
by juxtaposition, and $[a,b] = ab - ba$). Algebras satisfying both these identities are 
dubbed \textit{two-sided Alia algebras} in \cite{alia-d} and \cite{alia-db}.
Note that the class of two-sided Alia algebras contains both Lie algebras and
Novikov algebras (for the latter, see the end of \S \ref{simple}),
so it appears to be a natural and interesting class of algebras to study. 

Moreover, it is easy to see that an algebra $A$
is two-sided Alia if and only if the associated ``minus'' algebra $A^{(-)}$ with 
multiplication defined by the bracket $\liebrack$, is a Lie algebra
(in other words, $A$ is Lie-admissible), and multiplication
in $A$ can be written as 
\begin{equation}\label{plus}
ab = [a,b] + \varphi(a,b) ,\footnote{
Added May 5, 2017: The correct formula is: $ab = \frac 12 [a,b] + \varphi(a,b)$.
}
\end{equation}
where $\varphi$ is an
$A^{(-)}$-valued commutative $2$-cocycle on $A^{(-)}$.

Note, however, that also any commutative (nonassociative) algebra is two-sided Alia,
so the question of description of simple algebras in this class does not make much sense
without imposing additional conditions. One such natural condition is, in a sense,
opposite to the condition of commutativity of $A$ -- namely,
that the Lie algebra $A^{(-)}$ is simple. In this way we arrive to the problem of 
description of commutative $2$-cocycles on simple Lie algebras.

Second, commutative $2$-cocycles appear naturally in the description
of the second cohomology of current Lie algebras (\cite[Theorem 1]{without-unit}).
All this, as well as a sheer curiosity in what happens with the usual second Lie algebra
cohomology when we replace the condition of the skew-symmetricity of cochains by its opposite --
symmetricity -- makes them worth to study.

It is worth to note that this situation is similar (and somewhat dual) 
to the question which goes back to A.A. Albert and was a subject of an intensive
study later, namely, determination of Lie-admissible third power-associative algebras $A$
whose ``minus'' algebra $A^{(-)}$ belongs to some distinguished class of Lie algebras,
such as finite-dimensional simple Lie algebras, or Kac-Moody algebras.
It is easy to see that in terms of decomposition (\ref{plus}), 
for a given Lie algebra $L = A^{(-)}$ the third power-associativity implies the condition
$$
[\varphi(x,y),z] + [\varphi(z,x),y] + [\varphi(y,z),x] = 0 
$$
for any $x,y,z\in L$ (see, for example, \cite[p. 39]{benkart}).
The latter condition, together with (\ref{cocycle}), could be viewed as two parts 
of the usual $2$-cocycle equation on a Lie algebra with coefficients in the adjoint module.

Of course, the very definition of commutative $2$-cocycles begs for a proper generalization --
to define higher cohomology groups such that commutative $2$-cocycles constitute
cohomology of low degree. However, all ``naive'' attempts to construct such higher cohomology
seemingly fail. Maybe it could be developed in the framework of operadic cohomology of
two-sided Alia algebras.

The contents of the paper are as follows. In \S 1 we exhibit an exact sequence relating
commutative $2$-cocycles and so-called antiderivations, similar
to the well known relationship between the second cohomology and derivations with values
in the coadjoint module. 
In \S 2 we prove that in the class of simple Lie algebras
the property to possess nonzero commutative $2$-cocycles is equivalent to the property to 
satisfy the standard identity of degree $5$. In the class of 
finite-dimensional simple Lie algebras, this provides yet another characterization
of $sl(2)$ and the modular Zassenhaus algebra.
We also compute the space of commutative $2$-cocycles on the Zassenhaus algebra,
and note how this may be utilized in classification of simple Novikov algebras.
\S 3 is devoted to some speculations about the characteristic $3$ case. 
In \S 4 we establish formula expressing the space of commutative $2$-cocycles
on current Lie algebras in terms of its tensor factors, similar to the known formula
for the second cohomology. In the next two sections that formula is applied
to compute the space of commutative $2$-cocycles for various classes of
Lie algebras: Kac-Moody algebras in \S 5, and modular semisimple Lie algebras in \S 6.

Not surprisingly, the study of commutative $2$-cocycles turns out to be similar in some
aspects to the study of the second cohomology of Lie algebras, with some results having 
direct counterparts. There are, however, also significant differences: for example,
in the class of simple Lie algebras and close to them, 
non-trivial commutative $2$-cocycles turns out to be a very rare phenomenon.

\section*{Notation and conventions}

Throughout the paper, all algebras and vector spaces are defined over a ground field
$K$ of characteristic $\ne 2,3$, except in \S \ref{p3}.

For a given Lie algebra $L$, the space of all commutative $2$-cocycles on it will be
denoted as $Z^2_{comm}(L)$.
Occasionally, we will consider the \textit{$M$-valued commutative $2$-cocycles},
that is, bilinear forms $L \times L \to M$, where $M$ is an $L$-module (or, rather, just a 
vector space), satisfying the 
cocycle equation (\ref{cocycle}). It is immediate that if either $L$ or $M$ is 
finite-dimensional, then the space of all such cocycles
is isomorphic to $Z^2_{comm}(L) \otimes M$.

Recall that a Lie algebra $L$ is called \textit{perfect} if $[L,L] = L$.
Obviously, a symmetric bilinear form on $L$ that vanishes whenever one of the arguments belongs 
to $[L,L]$, is a commutative $2$-cocycle. Such cocycles will be called \textit{trivial}
and they exist on any non-perfect Lie algebra.
The space of trivial commutative $2$-cocycles is isomorphic to $S^2(L/[L,L])^*$.

The rest of our notation and definitions is mostly standard.

$H^n(L,M)$ and $Z^n(L,M)$ denote the usual Chevalley--Eilenberg $n$th order cohomology 
and the space of $n$th order cocycles, respectively, of a Lie algebra $L$ with coefficients in a module 
$M$.

When being considered as an $L$-module, the ground field $K$ is always understood as 
the one-di\-men\-si\-o\-nal trivial module.

$HC^1(A)$ denotes the first order cyclic cohomology of an associative commutative algebra
$A$. Recall that it is nothing but the space of all skew symmetric bilinear forms
$\alpha: A \times A \to K$ such that 
$$
\alpha(ab,c) + \alpha(ca,b) + \alpha(bc,a) = 0
$$
for any $a,b,c\in A$. Note an obvious but useful fact: if $A$ contains a unit $1$,
then $\alpha(1,A) = 0$ for any $\alpha\in HC^1(A)$.

$Der(A)$ denotes the Lie algebra of derivations of $A$.

For an associative algebra $A$, we may consider associated Lie $\liebrack$ and Jordan 
$\circ$ products on $A$: $[a,b] = ab - ba$ and $a \circ b = ab + ba$ for $a,b\in A$. 
The vector space $A$ with the bracket $\liebrack$ forms a Lie algebra which is denoted as 
$A^{(-)}$. 

For a vector space $V$, $S^2(V)$ and $\wedge^2(V)$ denotes the symmetric and skew-symmetric
product respectively (so, the space of symmetric, respectively skew-symmetric  bilinear 
forms on $V$ is isomorphic to $S^2(V)^*$, respectively to $\wedge^2(V)^*$).

Other nonstandard notions (antiderivations, cyclic forms) are introduced in \S 1.

\section{An exact sequence connecting commutative $2$-cocycles and antiderivations}

A relationship between the second cohomology with coefficients in the trivial module $H^2(L,K)$,
and the first cohomology with coefficients in the coadjoint module $H^1(L,L^*)$ 
was noted many times in slightly different forms in the literature, and goes back to the
classical works of Koszul and Hochschild--Serre. Namely,
there is an exact sequence
\begin{equation}\label{seq}
0 \to H^2(L,K) \overset{u}\to H^1(L,L^*) \overset{v}\to B(L)
\overset{w}\to H^3(L,K) .
\end{equation}
Here $B(L)$ denotes the space of symmetric bilinear invariant forms on $L$. 
The maps are defined as follows:
for the representative $\varphi\in Z^2(L,K)$ of a given
cohomology class, we have to take the class of $u(\varphi)$, the
latter being given by
\begin{equation}\label{u}
(u(\varphi)(x))(y) = \varphi(x,y)
\end{equation}
for any $x,y\in L$, $v$ is sending the class of a given cocycle $d\in Z^1(L,L^*)$ to the
bilinear form $v(d): L\times L\to K$ defined by the formula
\begin{equation*}
v(d)(x,y) = d(x)(y) + d(y)(x),
\end{equation*}
and $w$ is sending a given symmetric
bilinear invariant form $\varphi: L\times L \to K$ to the class of
the cocycle $\omega\in Z^3(L,K)$ defined by
\begin{equation}\label{w}
\omega(x,y,z) = \varphi([x,y],z)
\end{equation}
(see, for example, \cite{dzhu}, where a certain long exact sequence is obtained, of
which this one is the beginning, and references therein for many
earlier particular variations; this exact sequence was also
established in \cite[Proposition 7.2]{neeb} with two additional
terms on the right).

The following is a ``commutative'' version of this exact sequence, connecting
commutative $2$-cocycles and antiderivations. 

\begin{definition}
An \textit{antiderivation} of a Lie algebra
$L$ to an $L$-module $M$ is a linear map $d: L \to M$ such that 
$$
d([x,y]) = y \bullet d(x) - x \bullet d(y)
$$ 
for any $x,y \in L$, where $\bullet$ denotes
the $L$-action on $M$. The set of all such maps will be denoted as $ADer(L,M)$.
\end{definition}
 
When $M=L$, the adjoint module, we get the notion of an antiderivation 
of a Lie algebra (to itself) with the defining condition 
$$
d([x,y]) = - [d(x),y] - [x,d(y)] ,
$$ 
what was the subject of study in a number of papers, including \cite{filippov}. 

The third ingredient in our exact sequence, a counterpart of symmetric bilinear 
invariant forms, is defined as follows.

\begin{definition}
A bilinear form $\varphi:L \times L \to K$ is said to be
\textit{cyclic} if 
\begin{equation}\label{cyclic}
\varphi([x,y],z) = \varphi([z,x],y)
\end{equation}
for any $x,y,z\in L$. The space of all cyclic skew-symmetric bilinear forms on a Lie
algebra $L$ will be denoted as $C(L)$.
\end{definition}

\begin{proposition}\label{prop-seq}
For any Lie algebra $L$, there is an exact sequence
$$
0 \to Z^2_{comm}(L) \overset{u}\to ADer(L,L^*) \overset{v}\to C(L)
$$
where the map $u$ is defined by formula {\rm (\ref{u})}, and $v$ is sending a 
given antiderivation $d\in ADer(L,L^*)$ to the
bilinear form $v(d): L\times L\to K$ defined by the formula
\begin{equation*}
v(d)(x,y) = d(x)(y) - d(y)(x) .
\end{equation*}
\end{proposition}

\begin{proof}
If $\varphi \in Z^2_{comm}(L)$, then the cocycle equation yields 
\begin{multline*}
(u(\varphi)([x,y]))(z) = \varphi([x,y],z) = 
- \varphi(x,[y,z]) + \varphi(y,[x,z]) = 
(y \bullet u(\varphi)(x))(z) - (x \bullet u(\varphi)(y))(z)
\end{multline*}
for any $x,y,z \in L$, 
where $\bullet$ denotes the standard $L$-action on $L^*$, 
i.e. $u(\varphi) \in ADer (L, L^*)$.

If $d\in ADer(L,L^*)$, then $v(d)$ is obviously skew-symmetric, and 
\begin{multline*}
v(d)([x,y],z) = d([x,y])(z) - d(z)([x,y]) \\ = 
- d(x)([y,z]) + d(y)([x,z]) - d(z)([x,y]) = 
d([z,x])(y) - d(y)([z,x]) \\ = 
v(d)([z,x],y) 
\end{multline*}
for any $x,y,z \in L$, i.e. $v(d) \in C(L)$.

Now let us check exactness. Obviously, $u$ is injective, so 
the sequence is exact at $Z^2_{comm}(L)$. 
Next, $Im\,u$ consists of all antiderivations $d\in ADer(L,L^*)$ such that
the bilinear form $(x,y) \mapsto d(x)(y)$ satisfies the cocycle equation, 
what is equivalent to $d$ being antiderivation, and is symmetric, what is equivalent
to $d(x)(y) = d(y)(x)$ for any $x,y\in L$. But the latter condition is equivalent
to $d$ belonging to $Ker\,v$, hence $Im\,u = Ker\,v$ and the sequence is exact at 
$ADer(L,L^*)$.
\end{proof}

Unfortunately, we do not know how to extend this exact sequence further.
The map defined by the formula (\ref{w}) maps $C(L)$ to the space of 
skew-symmetric trilinear forms $L \times L \times L \to K$, but the resulting images 
are neither
Chevalley--Eilenberg (or Leibniz) cocycles, nor do they satisfy any other natural condition.
Obviously, this is related to the difficulty to define higher analogs
of commutative $2$-cocycles.

Another difficulty is related to the fact that $C(L)$ turns out to be not a very fascinating
invariant of a Lie algebra: it vanishes in the most interesting cases.

\begin{lemma}[Referee]\label{no-cl}
If $L$ is a perfect Lie algebra, then $C(L) = 0$.
\end{lemma}

\begin{proof}
For any $\varphi\in C(L)$, and any $x,y,z,t\in L$, we have
\begin{align*}
  \varphi([[x,y],z],t) &                           \\
\text{\small (Jacobi identity)} &= - \varphi([[z,x],y],t) - \varphi([[y,z],x],t)   \\
\text{\small (cyclicity of $\varphi$)} &= - \varphi([y,t],[z,x]) - \varphi([x,t],[y,z])   \\
\text{\small (skew-symmetry of $\varphi$)} &= - \varphi([x,z],[y,t]) + \varphi([y,z],[x,t])   \\
\text{\small (cyclicity of $\varphi$)} &= - \varphi([[y,t],x],z) + \varphi([[x,t],y],z)   \\
\text{\small (Jacobi identity)} &= \varphi([[x,y],t],z)                            \\
\text{\small (cyclicity of $\varphi$)} &= \varphi([z,[x,y]],t)                            \\
&= - \varphi([[x,y],z],t).
\end{align*}

Therefore, $\varphi([[L,L],L],L) = 0$, and the asserted statement follows.
\end{proof}

The exact sequence (\ref{seq}) was utilized many times in the literature
to evaluate $H^2(L,K)$ basing on $H^1(L,L^*)$ (see, for example, references in 
\cite{dzhu}). 
We will utilize Proposition \ref{prop-seq} in a similar way. For this, we shall
need to establish some facts about antiderivations of Lie algebras, 
what is done in the next section.

\section{Simple Lie algebras}\label{simple}

In \cite{filippov}, Filippov obtained many results about Lie algebras possessing 
a nonzero antiderivation. A slight modification of his reasonings allows to extend
some of these results to antiderivations with values in a Lie algebra module.

Recall the \textit{standard identity of degree $5$} in the class of Lie algebras:
\begin{equation}\label{deg5}
\sum_{\sigma\in S_4} 
(-1)^\sigma [[[[y, x_{\sigma(1)}], x_{\sigma(2)}], x_{\sigma(3)}], x_{\sigma(4)}]
= 0 ,
\end{equation}
where the summation is performed over all elements of the symmetric group $S_4$.
The word at the left side of (\ref{deg5}) will be denoted as $s_4(x_1,x_2,x_3,x_4,y)$.
For a Lie algebra $L$, let $s_4(L)$ denote the linear span of values of 
this word for any $x_1,x_2,x_3,x_4,y \in L$.
By the general result about verbal ideals (see, for example, \cite[Chapter 14]{as-book}, 
Theorem 2.8 and remark after it), $s_4(L)$ is an ideal of $L$.

\begin{lemma}\label{hom-image}
If a Lie algebra possesses nonzero commutative $2$-cocycles, then it has a 
nonzero homomorphic image satisfying the standard identity of degree $5$.
\end{lemma}

\begin{proof}
If a Lie algebra $L$ possesses nonzero commutative $2$-cocycles, then 
by Proposition \ref{prop-seq}, $L$ possesses a nonzero antiderivation to its 
coadjoint module. 
By \cite[Lemma 4.4]{delta} (which is just a slight generalization of 
\cite[Theorem 4]{filippov}), $s_4(L)$ lies in the kernel of this antiderivation,
and hence is a proper ideal in $L$. The required homomorphic image is $L/s_4(L)$.
\end{proof}

\begin{lemma}\label{codim1}
A Lie algebra with a subalgebra of codimension $1$ which is not an ideal, possesses
nonzero commutative $2$-cocycles.
\end{lemma}

\begin{proof}
Let $S$ be a subalgebra of codimension $1$ in a Lie algebra $L$. Pick 
$x \in L \backslash S$ and let $f: S \to K$ be a linear map such that
$[x,y] \in f(y)x + S$ for any $y\in S$. The Jacobi identity
implies $f([S,S]) = 0$. Since $S$ is not an ideal in $L$, $f$ is nonzero.
It is easy to see that any bilinear map $\varphi: L \times L \to K$ satisfying the
conditions
$$
\varphi(S,S) = 0, \quad \varphi(x,y) = \varphi(y,x) = f(y)
$$
for $y\in S$, is a commutative $2$-cocycle
(in fact, the cocycle equation in that case is equivalent to the Jacobi identity;
$\varphi(x,x)$ can take any value).
\end{proof}

\begin{theorem}\label{th-deg5}
A simple Lie algebra possesses nonzero commutative $2$-cocycles if and only if
it satisfies the standard identity of degree $5$.
\end{theorem}

\begin{proof}
The ``only if'' part follows from Lemma \ref{hom-image}.

The ``if'' part.
Suppose a simple Lie algebra $L$ satisfies the identity of degree $5$.
By Razmyslov's characterization of such algebras 
(\cite[Proposition 46.1]{razmyslov}), there exists a field extension $F$ of the centroid $C$
of $L$ such that the Lie $F$-algebra $L \otimes_C F$ contains a subalgebra of 
codimension $1$, and, by Lemma \ref{codim1}, possesses nonzero commutative
$2$-cocycles. As the space of commutative $2$-cocycles is obviously preserved under the
ground field extension, the Lie $C$-algebra $L \otimes_K C$ possesses nonzero commutative
$2$-cocycles. Each such cocycle, being restricted to $L$,
gives rise to a nonzero bilinear $K$-map
$L \times L \to C$, which is a commutative $C$-valued $2$-cocycle on $L$.
\end{proof}

Using other characterizations of simple Lie algebras satisfying the
standard identity of degree $5$, it is possible to give
an alternative proof of the ``if'' part of  Theorem \ref{th-deg5}
in the case where the ground field has positive characteristic.
Namely, according to \cite[Theorem 46.2]{razmyslov},
each simple Lie algebra of dimension $>3$ over the field of positive characteristic, 
satisfying the standard identity of degree $5$, can be represented
as a derivation algebra $A\partial = \set{a\partial}{a\in A}$ of an associative 
commutative algebra $A$ with unit, generated as an $A$-module by a single derivation 
$\partial\in Der(A)$. In addition, $A$ does not have $\partial$-invariant ideals.
For such algebras, we have

\begin{lemma}\label{ad}
$A^*$ is embedded into $Z^2_{comm}(A\partial)$.
\end{lemma}

\begin{proof}
Consider the map $F: A \to A\partial$ defined by $a\mapsto a\partial$ 
for $a\in A$. The kernel of this map is a $\partial$-invariant ideal of $A$, 
and hence this map is injective.

For any $\chi\in A^*$, define the map $\Phi(\chi): A\partial \times A\partial \to K$ as
$\Phi(\chi)(a\partial,b\partial) = \chi(ab)$ for $a,b\in A$. 
This map is well defined due to injectivity of $F$ and is a commutative $2$-cocycle:
\begin{multline*}
\Phi(\chi)([a\partial,b\partial],c\partial) + \Phi(\chi)([c\partial,a\partial],b\partial) + \Phi(\chi)([b\partial,c\partial],a\partial) \\ = 
\Phi(\chi)((a\partial(b) - b\partial(a))\partial,c\partial) \\+ \Phi(\chi)((c\partial(a) - a\partial(c))\partial,b\partial) \\ + \Phi(\chi)((b\partial(c) - c\partial(b))\partial,a\partial) \\ =
\chi(ac\partial(b)) - \chi(bc\partial(a)) + \chi(bc\partial(a)) - \chi(ab\partial(c)) + \chi(ab\partial(c)) - \chi(ac\partial(b)) = 0
\end{multline*}
for any $a,b,c\in A$. It is obvious that $\Phi$ is injective.
\end{proof}

\begin{corollary}\label{charact}
A finite-dimensional central simple Lie algebra possesses
nonzero commutative $2$-co\-cyc\-les if and only if it is isomorphic either to a form of $sl(2)$
or, in the case where the characteristic of the ground field is positive, to the 
Zassenhaus algebra $W_1(n)$.
\end{corollary}

\begin{proof}
Obviously, we may pass to the algebraic closure of the ground field.
Finite-dimensional simple Lie algebras over an algebraically closed field satisfying
the standard identity of degree $5$, are precisely $sl(2)$ and the Zassenhaus algebra. 
This well known fact
can be derived in several ways, perhaps the easiest one is to invoke once again 
\cite[Proposition 46.1]{razmyslov} to establish that such Lie algebras have a subalgebra
of codimension $1$, and then refer to the known classification of such algebras
(see, for example, \cite{elduque} and references therein).
\end{proof}

This generalizes \cite[Theorem 1.1]{alia-db}, where the same result is proved for Lie 
algebras of classical type by performing computations with the corresponding root system.
Another proof for such Lie algebras could be derived by combining results of 
\cite{larsson} and \cite{without-unit}. In a sense, root space computations
in \cite{alia-db} are equivalent to the appropriate part of computations in 
\cite{larsson}.

Note that the same reasoning (Theorem \ref{th-deg5} coupled with Razmyslov's
results) shows that among infinite-dimensional Lie 
algebras of Cartan type, only the Witt algebra possesses nonzero commutative $2$-cocycles.
The latters are described in \cite[Theorem 6.7]{alia-d}.
Another important class of infinite-dimensional Lie algebras -- 
Kac-Moody algebras -- is treated in \S \ref{kac-moody}.

The three-dimensional algebra $sl(2)$ and the Zassenhaus algebra are characterized among simple finite-di\-men\-si\-o\-nal
Lie algebras in various interesting ways: these are algebras having a subalgebra of 
codimension $1$, 
algebras having a solvable maximal subalgebra (see \cite[Corollary 9.1.0]{strade}), 
algebras with certain properties of the lattice of subalgebras 
(see, for example, \cite{varea} and references therein), algebras having
non-trivial $\delta$-derivations (\cite[\S 2]{delta}), etc.
Corollary \ref{charact} adds another characterization to this list.

Of course, it is interesting to compute exactly the space of commutative $2$-cocycles
on these algebras.
As noted in \cite[Theorem 6.5]{alia-d}, $Z^2_{comm}(sl(2))$ forms a $5$-dimensional
hyperplane in the $6$-dimensional space of symmetric bilinear maps 
$\varphi: sl(2) \times sl(2) \to K$, determined by the equation 
$$
\varphi(h,h) = 2\varphi(e_-,e_+)
$$
in the $sl(2)$-basis $\{e_-,e_+,h\}$ with multiplication table 
$[h,e_\pm] = \pm e_\pm$, $[e_-,e_+] = h$.

The famous Zassenhaus algebra $W_1(n)$ can be realized in two different ways. 
One realization 
is the algebra of derivations of the divided powers algebra 
$O_1(n) = \set{x^i}{0 \le i \le p^{n-1}}$, 
where $p$ is the characteristic of the ground field, with multiplication given by 
$$
x^i x^j = \binom{i+j}{j} x^{i+j} ,
$$
of the form
$O_1(n) \partial$, where the derivation $\partial \in Der(O_1(n))$ acts as follows:
$\partial (x^i) = x^{i-1}$.
Another realization is the algebra with the basis $\set{e_\alpha}{\alpha \in G}$,
where $G$ is an additive subgroup of order $p^n$ of the ground field $K$, with 
multiplication 
\begin{equation}\label{alphabeta}
[e_\alpha, e_\beta] = (\beta - \alpha) e_{\alpha + \beta}
\end{equation}
for $\alpha, \beta \in G$. We formulate the result in terms of the first realization,
but perform actual computations with the second one.

\begin{proposition}\label{zass}
$Z^2_{comm}(W_1(n)) \simeq O_1(n)^*$, each cocycle being of the form 
$$
(a\partial, b\partial) \mapsto f(ab)
$$ 
where $a,b \in O_1(n)$, for some $f\in O_1(n)^*$.
\end{proposition}

\begin{proof}
The proof consists of straightforward computations reminiscent both the computation of the
second cohomology of $W_1(n)$ in \cite[\S 5]{block} and computation of the space
of commutative $2$-cocycles on the infinite-dimensional Witt algebra in 
\cite[Theorem 6.7]{alia-d}.

Let $\varphi \in Z^2_{comm}(W_1(n))$. Writing the cocycle equation in terms of the basis
with multiplication (\ref{alphabeta}) for 
$e_0, e_\alpha, e_\beta$, where $\alpha, \beta \in G$ such that $\alpha \ne \beta$, we get
\begin{equation}\label{yoyo}
\varphi(e_\alpha, e_\beta) = \varphi(e_0, e_{\alpha + \beta}) .
\end{equation}
Now, taking this into 
account, and writing the cocycle equation for $e_\alpha, e_\beta, e_{\alpha + \beta}$,
where $\alpha, \beta \in G$ such that $\alpha, \beta \ne 0$, we get
$$
\varphi (e_{\alpha + \beta}, e_{\alpha + \beta}) = \varphi (e_0, e_{2\alpha + 2\beta}) .
$$
Consequently, (\ref{yoyo}) holds for any $\alpha, \beta \in G$. Conversely, each 
symmetric bilinear map satisfying the condition (\ref{yoyo}) is easily seen to satisfy
the cocycle equation. Each such map can be decomposed into the sum of the maps
of the form 
$$
(e_\alpha, e_\beta) \mapsto \begin{cases}
1, & \alpha + \beta = \gamma \\
0, & \text{otherwise}
\end{cases}
$$
for each $\gamma \in G$. Thus we get $|G| = p^n$ linearly independent cocycles, so
the whole space $Z^2_{comm}(W_1(n))$ is $p^n$-dimensional. Now, switching to the first realization as derivation 
algebra of $O_1(n)$, applying Lemma \ref{ad} and comparing dimensions, we see
that the embedding of Lemma \ref{ad} is an isomorphism in this particular case.
\end{proof}

Note that the results of this section allow to provide a somewhat alternative way 
for classification of finite-dimensional simple Novikov algebras. 
Recall that an algebra is called
\textit{Novikov} if it satisfies the identities
$$
x(yz) - (xy)z = y(xz) - (yx)z
$$
and
$$
(xy)z = (xz)y .
$$
Novikov algebras are ubiquitous in various branches of mathematics and physics
(see \cite{burde}, \cite{osborn} and \cite{zelmanov} with a transitive closure of 
references therein)\footnote[2]{
Actually, what we define here are \textit{left Novikov} algebras. 
\textit{Right Novikov} algebras are defined by the opposite identities,
interchanging left and right multiplications.
Often in the literature,
left Novikov and right Novikov are confused, even if authors make an explicit
attempt not to do so (for example, \cite{zelmanov} treats right Novikov algebras, while
\cite{osborn} and \cite{burde} -- left Novikov ones). 
However, almost every result about left Novikov algebras is easily
transformed into one about right Novikov algebras, and vice versa.
}.

A well known important fact is that each Novikov algebra $A$ is Lie-admissible, i.e.,
$A^{(-)}$ is a Lie algebra.
In \cite{zelmanov}, Zelmanov proved that finite-dimensional simple Novikov algebras over a field
of characteristic zero are $1$-dimensional (i.e., coincide with the ground field),
and in \cite{osborn}, Osborn proved that if $A$ is a finite-dimensional simple Novikov 
algebra over a field of positive characteristic, then $A^{(-)}$ is either 
$1$-dimensional, or isomorphic to the Zassenhaus algebra. 
The latter was the key result in obtaining later a complete classification of simple 
Novikov algebras.

Another important observation -- which is a matter of simple calculations (see, for example, 
\cite[Lemma 2.3]{burde}) -- is that Novikov algebras are two-sided Alia.
Further, it follows from the proof of \cite[Theorem 3.5]{osborn}, that if $A$ is a 
finite-dimensional simple Novikov algebra, then the Lie algebra $A^{(-)}$ is also 
simple. Thus, as noted in the Introduction, the multiplication in a finite-dimensional 
simple Novikov algebra $A$ could be written in the form (\ref{plus}), where $\liebrack$
defines a simple Lie algebra structure on $A$, and $\varphi$ defines a commutative
$2$-cocycle on that Lie algebra.  Note further that $\varphi$ is a nonzero cocycle,
as otherwise $A$ is a Lie algebra, what quickly leads to contradiction.

Then, instead of appealing to some structural results about graded and filtered Lie
algebras, like in \cite{osborn}, we may invoke Corollary \ref{charact} to see that
$A^{(-)}$ is isomorphic either to $1$-dimensional abelian algebra, or to $sl(2)$, or to 
the Zassenhaus algebra. 
Straightforward calculations based on \cite[Theorem 6.5]{alia-d} exclude the case of 
$sl(2)$, and in this way we get the main results of \cite{zelmanov} and \cite{osborn}.
Proposition \ref{zass} may be utilized in the subsequent classification.

\section{Characteristic $3$}\label{p3}

The case of characteristic $3$ differs drastically:
as noted in \cite[\S 6.2]{alia-d}, in this case
any symmetric bilinear invariant form on a Lie algebra is a commutative $2$-cocycle.
Let us rephrase this observation in another trivial, yet not without interest, form.

To start with, let us make a useful observation valid in any characteristic.
For a Lie algebra $L$, consider the standard $Der(L)$-action on the space of 
bilinear forms on $L$: for $\varphi$ being such a form, and $D\in Der(L)$,
\begin{equation}\label{act}
(D \bullet \varphi)(x,y) = -\, \varphi(D(x),y) - \varphi(x,D(y)) ,
\end{equation}
$x,y\in L$.
It is well known (and easy to verify) that the space of symmetric bilinear invariant forms on 
$L$ is closed under this action. 
It turns out that the same is true for the space of commutative $2$-cocycles:

\begin{lemma}\label{invar}
For any Lie algebra $L$, $Z^2_{comm}(L)$ is closed under the action (\ref{act}).
\end{lemma}

\begin{proof}
This follows from the following equality, valid for any bilinear form 
$\varphi: L \times L \to K$, any $D\in Der(L)$, and any $x,y,z\in L$:
$$
d(D\bullet \varphi)(x,y,z) = 
-\, d\varphi(x,y,D(z)) 
- d\varphi(z,x,D(y)) 
- d\varphi(y,z,D(x)) .
$$
where $d\varphi(x,y,z)$ denotes the left-hand side of equality (\ref{cocycle}).
\end{proof}

In particular, $Z^2_{comm}(L)$ is closed under the action of $L$ (via inner derivations).
This is, essentially, the same observation which is used in deriving a very useful fact 
about triviality of the Lie algebra action on its cohomology. Note, however, that,
unlike for cohomology, in the case of commutative $2$-cocycles
we do not have coboundaries, so we cannot normalize 
the cocycle appropriately to derive the triviality of this action.
Moreover, for invariants of this action we have the following obvious dichotomy:
$$
Z^2_{comm}(L)^L = \begin{cases}
B(L),                                & p=3    \\
\set{\varphi: L \times L \to K}{\varphi \text{ is symmetric, } \varphi([L,L],L) = 0}, & p\ne 3 .
\end{cases}
$$

Still, we can make use of Lemma \ref{invar} in the same way as in cohomological
considerations: when $T$ is a torus in a Lie algebra $L$ such that $L$ decomposes into
the direct sum of the root spaces with respect to the action of $T$
(what always takes places if $L$ is finite-dimensional and the ground field
is algebraically closed), then $Z^2_{comm}(L)$ decomposes into the direct sum
of root spaces with respect to the induced $T$-action. 

Naturally, in order to compute the space of commutative $2$-cocycles on
any class of Lie algebras of characteristic $3$, it would be beneficial
to elucidate first what the space of symmetric bilinear invariant forms on these algebras looks like. 

\begin{conjecture}\label{conj-symm}
The space of symmetric bilinear invariant forms on any central Lie algebra of classical 
type over a field of arbitrary characteristic is $1$-dimensional.
\end{conjecture}

Note that in small characteristics (including characteristic $3$) not all Lie algebras obtained via the usual 
Chevalley basis construction, are simple (see, for example, \cite[\S 4.4]{strade}). 
In such cases, under Lie algebras of classical type we mean both these 
non-simple Lie algebras, and their central simple quotients.

Of course, for $p=0$ and $p>5$ this conjecture trivially follows from the well 
known statement about the Killing form. However, in small characteristics the Killing form, and, more general,
any trace form, may vanish (see \cite{garibaldi} and references therein), 
so this conjecture, perhaps, comes as a bit of surprise. It is, however, supported
by computer calculations\footnote[2]{
A simple-minded GAP code for calculation of the spaces of 
commutative $2$-co\-cyc\-les and of symmetric bilinear invariant forms on a Lie
algebra, is available at \texttt{http://justpasha.org/math/comm2.gap} . 
[Added January 4, 2016: currently available as an ancillary file accompanying 
the arXiv version of this paper]. The code works by writing the corresponding conditions in terms of a certain basis of a Lie algebra, 
and solving the arising linear homogeneous system.
}.

\begin{conjecture}\label{th-p3}
The space of commutative $2$-cocycles on a finite-dimensional Lie algebra of 
classical type different from $sl(2)$, over the field of characteristic $3$,
coincides with the space of symmetric bilinear invariant forms (and, hence,
by Conjecture \ref{conj-symm}, is $1$-dimensional).
\end{conjecture}

Note that for Lie algebras not of classical type, it is no longer true that 
every commutative $2$-cocycle is a symmetric bilinear invariant form: for example, 
for the Zassenhaus algebra $W_1(n)$, Proposition \ref{zass} is valid also in characteristic 
$3$ (with, essentially, the same proof). On the other hand, since $W_1(n)$ is simple,
the space of symmetric bilinear invariant forms on it is at most $1$-dimensional
(it is, in fact, $1$-dimensional in characteristic $3$, as shown in
\cite[\S 2, Corollary]{dzhu-izv} or \cite[\S4.6, Theorem 6.3]{sf}).

A straightforward way to prove Conjectures \ref{conj-symm} and \ref{th-p3} seems to 
employ the strategy of \cite[Theorem 1.1]{alia-db}: first to prove the statements
for algebras of rank $2$, and then derive the general case.

Lie algebras of classical type of rank $2$ are $A_2$, $B_2$ and $G_2$. For these algebras,
one may verify on a computer that both the spaces of commutative $2$-cocycles (for $p=3$)
and of symmetric bilinear invariant forms (for any $p$), are $1$-dimensional.
The general case should follow then from considerations of $2$-sections 
in a Cartan decomposition, with the help of Lemma \ref{invar}.

However, there are lot of subtleties when dealing with structure constants of classical
Lie algebras in small characteristics, as demonstrated by a noticeable amount of errors
in works devoted to such algebras. We postpone this laborious task to the future.

\section{Current Lie algebras}

In this section we consider the current Lie algebras, i.e. Lie algebras of the form
$L\otimes A$ where $L$ is a Lie algebra and $A$ is an associative commutative algebra,
with multiplication defined by 
$$[x\otimes a, y\otimes b] = [x,y]\otimes ab$$ 
for $x,y \in L$, $a,b \in A$.

\begin{theorem}\label{curr}
Let $L$ be a Lie algebra, $A$ an associative commutative algebra, 
and at least one of $L$, $A$ is finite-dimensional.
Then each commutative $2$-cocycle on $L\otimes A$ can be
represented as a sum of decomposable cocycles $\varphi\otimes \alpha$, 
$\varphi: L\times L \to K$, $\alpha: A\times A\to K$ of one of the $8$ following types:
\begin{enumerate}
\item 
$\varphi([x,y],z) + \varphi([z,x],y) + \varphi([y,z],x) = 0$ and 
$\alpha(ab,c) = \alpha(ca,b)$,
\item 
$\varphi([x,y],z) = \varphi([z,x],y)$ and $\alpha(ab,c) + \alpha(ca,b) + \alpha(bc,a) = 0$,
\item $\varphi ([L,L],L) = 0$,
\item $\alpha (AA,A) = 0$,
\end{enumerate}
and each of these four types splits into two subtypes: with both $\varphi$ and $\alpha$
symmetric, and with both $\varphi$ and $\alpha$ skew-symmetric.
\end{theorem}

\begin{proof}
The proof goes almost verbatim to the proof of Theorem 1 in \cite{without-unit}.
\end{proof}

\begin{corollary}\label{cor-curr}
Suppose all assumptions of Theorem \ref{curr} hold, and, additionally, 
$A$ contains a unit. Then
\begin{align*}
Z^2_{comm}(L\otimes A) \simeq\> & Z^2_{comm}(L) \otimes A^* 
\\ \oplus\> & C(L) \otimes HC^1(A)
\\ \oplus\> & (S^2(L/[L,L]))^* \otimes Ker(S^2(A) \to A)^* 
\\ \oplus\> & (\wedge^2(L/[L,L]))^* \otimes \set{ab \wedge c + ca \wedge b + bc \wedge a}{a,b,c\in A}^* ,
\end{align*}
where the map $S^2(A) \to A$ is induced by multiplication in $A$.
\end{corollary}

\begin{proof}
Consider the cases of Theorem \ref{curr}.

Case (i). 
Substituting $c=1$ in the condition $\alpha(ab,c) = \alpha(ca,b)$, we get
$\alpha(a,b) = \beta(ab)$ for a certain linear form $\beta: A \to K$.
Consequently, both $\varphi$ and $\alpha$ are necessarily symmetric,
hence $\varphi \in Z^2_{comm}(L)$ and the linear span of cocycles of this type is isomorphic
to $Z^2_{comm}(L) \otimes A^*$.

Case (ii).
Similarly, substituting $b=c=1$ in the condition 
$$
\alpha(ab,c) + \alpha(ca,b) + \alpha(bc,a) = 0 ,
$$
and assuming that $\alpha$ is symmetric, we get $\alpha(a,1) = 0$. Now substituting just $c=1$
in the same condition, we get that $\alpha$ vanishes. Consequently, for cocycles of 
this type both $\varphi$ and $\alpha$ are necessarily skew-symmetric, hence
$\varphi \in C(L)$ and $\alpha \in HC^1(A)$.

Case (iv). If $A$ contains a unit, then $AA = A$, so cocycles of this type vanish.

Singling out from the linear span of the remaining cocycles of type (iii) 
the direct sum complement to the linear span of cocycles of type (i) and (ii), 
and rearranging it in an obvious way, we get the desired isomorphism.
\end{proof}

Evidently, the last two direct summands at the right-hand side of the isomorphism of 
Corollary \ref{cor-curr} constitute trivial cocycles.

It is possible to extend Theorem \ref{curr} and Corollary \ref{cor-curr} to 
various generalizations of current Lie algebras, such as
twisted algebras, extended affine Lie algebras, toroidal Lie algebras, Lie algebras
graded by root systems, etc.
Some of computations could be quite cumbersome, but all of them seem to be
amenable to the technique used in \cite{without-unit}. 

It is possible also to consider a sort of noncommutative version of Corollary \ref{cor-curr}, 
namely, commutative $2$-cocycles on the Lie algebra $sl(n,A)$ for an associative 
(and not necessarily commutative) algebra $A$ with unit.
It is possible to show that any homomorphic image of such an algebra is closely related
to the algebra of the form $sl(n,B)$, where $B$ is a homomorphic image of $A$. 
If $sl(n,A)$ possesses nonzero 
$2$-commutative cocycles, then by Lemma \ref{hom-image}, at least one of these homomorphic 
images satisfies the standard identity of degree $5$. 
As $sl(n)$ is a subalgebra of $sl(n,B)$, this implies
that there are no commutative $2$-cocycles if $n>2$. The algebra $sl(2,A)$, on the contrary, possesses nonzero $2$-commutative cocycles, 
but in general there seems no nice expression for them in terms of $A$.
As the final answer turns out not to be very interesting in either case, 
we are not going into details.

\section{Kac-Moody algebras}\label{kac-moody}

If we want to apply the results of the preceding section to Kac-Moody algebras, 
we should deal not with the current Lie algebras and their twisted analogs, 
but their extensions by means of central elements and derivations.
To this end, we make the following elementary observations.

\begin{lemma}\label{6}
Let $L$ be a Lie algebra and $I$ an ideal of $L$. Then $Z^2_{comm}(L/I)$ is embedded
into $Z^2_{comm}(L)$.
\end{lemma}

\begin{proof}
There is an obvious bijection between $Z^2_{comm}(L/I)$ and the set of cocycles \newline
$\varphi\in Z^2_{comm}(L)$ such that $\varphi(L,I) = 0$.
\end{proof}

Of course, the similar embedding exists for ordinary (skew-symmetric) cocycles, but 
this embedding is usually not preserved on the level of cohomology. 
For example, free Lie algebras possess, in a sense, ``the most'' of $2$-cocycles,
accumulating all cocycles from their homomorphic images, but 
in the skew-symmetric case all cocycles are killed off by coboundaries. This
has no parallel in the commutative case due to absence of 
``commutative $2$-coboundaries''.

The following is a very particular complement, in a sense, to Lemma \ref{6}.

\begin{lemma}\label{perfect}
Let $L$ be a Lie algebra and $I$ a perfect ideal of codimension $1$ of $L$. 
Then $Z^2_{comm}(L)$ is embedded into $Z^2_{comm}(I) \oplus K$, the second direct summand
being represented by a trivial cocycle.
\end{lemma}

\begin{proof}
Write $L = I \oplus Kx$ for some $x\in L\backslash I$. Let $\varphi \in Z^2_{comm}(L)$.
The cocycle equation on $L$ is equivalent to the cocycle equation on $I$, plus the cocycle
equation for $a,b\in I$ and $x$, the latter could be written in the form
$$
\varphi([a,b],x) = \varphi([a,x],b) - \varphi([b,x],a) .
$$
As $[I,I] = I$, this formula provides a well defined extension of $\varphi$ from 
$I \times I$ to $I \times Kx$.

There are no further restrictions on the value of $\varphi(x,x)$, so defining a bilinear form 
$\varphi$ on $L$ by $\varphi(x,x) = 1$, $\varphi(I,I) = \varphi(I,x) = 0$, we get a
trivial cocycle which corresponds to the second direct summand.  
\end{proof}

We use the realization of affine Kac-Moody Lie algebras as extensions of non-twisted and twisted 
current Lie algebras (see \cite[Chapters 7 and 8]{kac}). Let 
$\mathfrak g = \bigoplus_{j\in \mathbb Z_n} \mathfrak g_j$ be a simple
finite-dimensional $\mathbb Z_n$-graded Lie algebra (it is possible that 
$\mathfrak g = \mathfrak g_0$, i.e., $n=1$ and the grading is trivial, what corresponds to the 
non-twisted case).
Consider the following subalgebra of the current Lie algebra $\mathfrak g \otimes K[t,t^{-1}]$:
$$
L(\mathfrak g, n) = \bigoplus_{i\in \mathbb Z} \mathfrak g_{i(mod\, n)} \otimes t^i .
$$
This algebra has a non-split perfect central extension $\widetilde{L}(\mathfrak g, n)$. 
Each affine Kac-Moody algebra can be represented in the form 
$$
\widehat{L}(\mathfrak g, n) = \widetilde{L}(\mathfrak g, n) \oplus Kt\frac{d}{dt}
$$
for a suitable $\mathfrak g$, where multiplication between elements
of $\widetilde{L}(\mathfrak g, n)$ and $t\frac{d}{dt}$ is defined by the action of the latter
as derivation on the algebra of Laurent polynomials $K[t,t^{-1}]$.

We will consider first the case of non-twisted affine Kac-Moody algebras in a bit more
general situation.
Let $L$ be a Lie algebra, 
$\langle \cdot, \cdot \rangle$ a nonzero symmetric bilinear invariant form on $L$,
$A$ a commutative associative algebra with unit,
$D$ a Lie subalgebra of $Der(A)$, and $\xi$ a nonzero $D$-invariant element of $HC^1(A)$.
Consider a Lie algebra defined as the vector space $(L \otimes A) \oplus Kz \oplus D$ 
with the following multiplication:
\begin{align*}
[x\otimes a, y\otimes b] &= [x,y] \otimes ab + \langle x, y \rangle \xi(a,b) z  \\
[x \otimes a, d] &= x \otimes d(a)
\end{align*}
for $x,y \in L$, $a,b \in A$, $d\in D$, and $z$ belongs to the center.
Note that the semidirect sum $(L\otimes A) \inplus D$ is a quotient by the $1$-dimensional central ideal $Kz$,
and $(L\otimes A) \oplus Kz$ is a subalgebra, the latter being central extension of the
current Lie algebra $L\otimes A$ (for generalities about central extensions of current Lie 
algebras, see \cite{asterisque, without-unit, neeb}).

Specializing this construction to the case where $K$ is an algebraically closed field of
characteristic zero, $L = \mathfrak g$, a simple finite-dimensional Lie algebra, 
$\langle \cdot,\cdot \rangle$ is the Killing form on $\mathfrak g$, $A = K[t,t^{-1}]$, 
$D = Kt\frac{d}{dt}$, and $\xi (f,g) = Res(g \frac{df}{dt})$ for $f,g\in K[t,t,^{-1}]$,
we get non-twisted affine Kac-Moody algebras.

\begin{lemma}\label{lemma-yaya}
Let $L$ be perfect, and one of $L$, $A$ is finite-dimensional. Then\footnote[2]{
Added May 24, 2011: By Lemma \ref{no-cl}, the term containing $C(L)$ is redundant here.
The proofs of Lemmata \ref{lemma-yaya} and \ref{lemma-yoyo} can be simplified 
accordingly. 
}
\begin{align*}
Z^2_{comm} &((L \otimes A) \inplus D) \\ \simeq\>
& Z^2_{comm}(L) \otimes \set{\chi\in A^*}{\chi(d(a)b - ad(b)) = 0 \text{ for any } a,b\in A, d\in D} 
\\ \oplus\>
& C(L) \otimes \set{\beta \in HC^1(A)}{\beta(d(a),b) - \beta(a,d(b)) = 0 
\text{ for any } a,b\in A, d\in D} 
\\ \oplus\> & Z^2_{comm}(D) .
\end{align*}
\end{lemma}

\begin{proof}
Let $\Phi \in Z^2_{comm} ((L \otimes A) \inplus D)$. A restriction of $\Phi$ to 
$(L\otimes A) \times (L \otimes A)$ is a commutative $2$-cocycle on $L\otimes A$. 
By Corollary \ref{cor-curr}, there are
$\varphi \in Z^2_{comm}(L)$, $\chi \in A^*$, $\alpha \in C(L)$ and $\beta \in HC^1(A)$
such that
\begin{equation}\label{phi}
\Phi (x\otimes a, y\otimes b) = \varphi(x,y) \chi(ab) + \alpha(x,y) \beta(a,b)
\end{equation}
for any $x,y \in L$, $a,b \in A$.
Writing the cocycle equation for $x\otimes a$, $y\otimes b$, $d$, we get
\begin{equation}\label{yaya}
\Phi ([x,y] \otimes ab, d) =   
\varphi(x,y) \chi(d(a)b - ad(b)) + \alpha(x,y) (\beta(d(a),b) - \beta(a,d(b)))
\end{equation}
for any $x,y\in L$, $a,b\in A$, and $d\in D$. Substituting here $a=1$ and $b=1$,
we get respectively:
\begin{align*}
\Phi([x,y] \otimes b, d) &= - \varphi(x,y) \chi(d(b)) \\
\Phi([x,y] \otimes a, d) &= \phantom{-}\varphi(x,y) \chi(d(a)) ,
\end{align*}
what implies $\Phi(L\otimes A, D) = 0$.
Hence the right-hand side of (\ref{yaya}) vanishes, and symmetrizing it further with respect 
to $x$, $y$, we get that both summands vanish separately.
Thus we see that $\chi\in A^*$ and $\beta\in HC^1(A)$ in formula (\ref{phi}) satisfy
additional conditions $\chi(d(a)b - ad(b)) = 0$ and $\beta(d(a),b) - \beta(a,d(b)) = 0$ 
for any $a,b\in A$.

The cocycle equation for one element from $L\otimes A$ and two elements from $D$ is 
satisfied trivially, and restriction of $\Phi$ to $D \times D$ gives rise to an
element from $Z^2_{comm}(D)$.
\end{proof}

\begin{lemma}\label{lemma-yoyo}
Let $L$ be perfect, and one of $L$, $A$ is finite-dimensional. Then
$$
Z^2_{comm} ((L \otimes A) \oplus Kz \inplus D) \simeq Z^2_{comm} ((L \otimes A) \inplus D) .
$$
\end{lemma}

\begin{proof}
By Lemma \ref{6}, we have an embedding of $Z^2_{comm} ((L \otimes A) \inplus D)$ 
into \newline $Z^2_{comm} ((L \otimes A) \oplus Kz \inplus D)$, and this embedding is an isomorphism 
if and only if 
\begin{equation}\label{iff}
\Phi((L \otimes A) \oplus Kz \inplus D, z) = 0
\end{equation}
for any $\Phi\in Z^2_{comm} ((L \otimes A) \oplus Kz \inplus D)$.

Writing the cocycle equation for $x\otimes a, y\otimes b, z$, we get
\begin{equation}\label{yoyo-1}
\Phi([x,y] \otimes ab + \langle x, y \rangle \xi(a,b)z, z) = 0
\end{equation}
for any $x,y\in L$, $a,b \in A$. 
Substituting here $b=1$, we get 
\begin{equation}\label{yoyo-2}
\Phi(L \otimes A, z) = 0 .
\end{equation}
Substituting the latter equality back to (\ref{yoyo-1}), we get $\Phi(z,z) = 0$.

Being restricted to $((L\otimes A) \oplus Kz) \times ((L\otimes A) \oplus Kz)$, $\Phi$
gives rise to a commutative $2$-cocycle on $(L\otimes A) \oplus Kz$, and due to 
(\ref{yoyo-2}), to a commutative $2$-cocycle on $L\otimes A$. Now we proceed as in 
the proof of Lemma \ref{lemma-yaya}: by Corollary \ref{cor-curr}, the equality
(\ref{phi}) holds for any $x,y\in L$, $a,b\in A$ and appropriate 
$\varphi\in Z^2_{comm}(L)$, $\chi\in A^*$, $\alpha\in C(L)$, and $\beta\in HC^1(A)$. 
Then the cocycle equation for $x\otimes a$, $y\otimes b$, $d$, gives
\begin{multline}\label{yaya1}
\Phi([x,y] \otimes ab, d) + \langle x,y \rangle \xi(a,b)\Phi(z,d) \\ = 
\varphi(x,y) \chi(d(a)b - ad(b)) + \alpha(x,y) (\beta(d(a),b) - \beta(a,d(b)))
\end{multline}
for any $x,y\in L$, $a,b\in A$, $d\in D$, and substitution of units in this equality 
yields $\Phi(L \otimes A, D) = 0$. Substituting this back to (\ref{yaya1}) and
symmetrizing with respect to $x,y$, gives
\begin{equation*}
\langle x,y \rangle \xi(a,b)\Phi(z,d) = \varphi(x,y) \chi(d(a)b - ad(b))
\end{equation*}
for any $x,y\in L$, $a,b\in A$, $d\in D$. Either both sides of this equality vanishes,
in which case $\Phi(z,d) = 0$ for any $d\in D$, or $\varphi(x,y) = \lambda \langle x,y \rangle$
for some nonzero $\lambda\in K$. But the latter is obviously impossible
(note that the condition that the characteristic of the ground field is different from $3$
is crucial here: as noted in \cite[\S 6.2]{alia-d}, in characteristic $3$ every
symmetric bilinear invariant form is a commutative $2$-cocycle).

Therefore, (\ref{iff}) holds and the desired isomorphism follows.
\end{proof}

Affine Kac-Moody algebras, being non-perfect, possess trivial nonzero commutative 
$2$-cocycles.
As the commutant is always of codimension $1$, the space of such cocycles is $1$-dimensional.
There are no other cocycles, as the following result shows.

\begin{theorem}
Affine Kac-Moody algebras do not possess non-trivial commutative $2$-cocycles.
\end{theorem}

\begin{proof}
First consider the case of non-twisted algebras.
By Lemmata \ref{lemma-yaya} and \ref{lemma-yoyo}, the space of commutative $2$-cocycles
on a non-twisted affine Kac-Moody algebra 
$$
(\mathfrak g \otimes K[t,t^{-1}]) \oplus Kz \inplus Kt\frac{d}{dt}
$$
is isomorphic to 
\begin{align*}
& Z^2_{comm}(\mathfrak g) \otimes 
\set{\chi\in K[t,t^{-1}]^*}{\chi(t\frac{df}{dt}g - ft\frac{dg}{dt}) = 0 \text{ for any } f,g\in K[t,t^{-1}]} 
\\ \oplus &
C(\mathfrak g) \otimes 
\set{\beta \in HC^1(K[t,t^{-1})}{\beta(t\frac{df}{dt},g) - \beta(f,t\frac{dg}{dt}) = 0 
\text{ for any } f,g\in K[t,t^{-1}]}
\\ \oplus &
Z^2_{comm}(Kt\frac{d}{dt}) .
\end{align*}

Let us look at the second tensor factor in the first direct summand.
Substituting $f=t^i$, $g=t^j$ in the defining condition 
\begin{equation*}
\chi(t\frac{df}{dt}g - ft\frac{dg}{dt}) = 0
\end{equation*}
for $\chi\in K[t,t^{-1}]^*$, we get $(i-j) \chi(t^{i+j}) = 0$ for any $i,j \in \mathbb Z$
such that $i \ne j$. Hence $\chi = 0$, and the first direct summand vanishes.

By Lemma \ref{no-cl}, $C(\mathfrak g) = 0$, hence the second direct summand 
vanishes too.

Finally, $Z^2_{comm}(Kt\frac{d}{dt})$, being the space of commutative cocycles on the 
$1$-di\-men\-si\-o\-nal Lie algebra, is $1$-di\-men\-si\-o\-nal and constitute trivial cocycles.

Now consider the general (twisted) case. 
Suppose that $\widehat{L}(\mathfrak g, n)$ possesses non-trivial commutative $2$-cocycles. Then by
Lemma \ref{perfect}, $\widetilde{L}(\mathfrak g, n)$ possesses nonzero commutative $2$-cocycles, and by Lemma \ref{hom-image}, $\widetilde{L}(\mathfrak g, n)$ has a 
nonzero homomorphic image satisfying the standard identity of degree $5$.
Every homomorphic image of $\widetilde{L}(\mathfrak g, n)$ is either a homomorphic image
of $L(\mathfrak g, n)$, or is a central extension of such.
In the latter case, factoring by the central element, we will get again a homomorphic
image of $L(\mathfrak g, n)$ satisfying the standard identity of degree $5$.
But $\mathfrak g$ is a homomorphic image of $L(\mathfrak g,n)$ under the evaluation homomorphism
$$
\sum_{i\in \mathbb Z} x_i \otimes t^i \mapsto \sum_{i\in \mathbb Z} x_i ,
$$
where $x_i\in \mathfrak g$, and all but a finite number of summands vanish.
Thus $\mathfrak g$ satisfies the standard identity of degree $5$ too, 
hence $\mathfrak g \simeq sl(2)$, and the corresponding Kac-Moody algebra is of type 
$A_1^{(1)}$, i.e. is a non-twisted algebra
of the form 
$$
(sl(2) \otimes K[t,t^{-1}]) \oplus Kz \inplus Kt\frac{d}{dt} .
$$
But the case of non-twisted algebras was already covered.
\end{proof}

Note that there is some ambiguity in definition of affine Kac-Moody algebras,
and sometimes they are defined without employing the derivation extension, i.e. merely as 
central extensions of non-twisted or twisted current algebras
(though some people argue that these are not ``real'' Kac-Moody algebras, as they
are not Lie algebras corresponding to Cartan matrices). 
Such algebras are perfect, and possess 
nonzero commutative $2$-cocycles only in the case of non-twisted type $A_1^{(1)}$.
The proof is absolutely similar to those presented above.

According to Lemma \ref{lemma-yoyo} (with $D = 0$) and Corollary \ref{cor-curr}, the space
of commutative $2$-cocycles on the Kac-Moody algebra $(sl(2) \otimes K[t,t^{-1}]) \oplus Kz$
of type $A_1^{(1)}$ is infinite-dimensional and is isomorphic to $Z^2_{comm}(sl(2)) \otimes K[t,t^{-1}]^*$, 
each cocycle being of the form
\begin{align*}
(x \otimes f, y \otimes g) & \mapsto \varphi(x,y) \otimes \chi(fg) \\
(x \otimes f, z) & \mapsto 0 \\
(z,z) & \mapsto 0
\end{align*}
where $x,y \in sl(2)$, $f,g \in K[t,t^{-1}]$, for some $\varphi \in Z^2_{comm}(sl(2))$,
$\chi \in K[t,t^{-1}]^*$.

\section{Modular semisimple Lie algebras} 

Essentially the same approach as in the previous section, allows to compute the space of 
commutative $2$-cocycles
on finite-dimensional semisimple Lie algebras over the field of positive characteristic $p$. 
According to Block's classical theorem 
(see, for example, \cite[Corollary 3.3.6]{strade}),
the typical examples of such algebras are Lie algebras of the form
$(S \otimes O_n) \inplus D$ where $S$ is a simple Lie algebra, 
$O_n = K[x_1, \dots, x_n]/(x_1^p, \dots, x_n^p)$ is the reduced polynomial algebra in $n$ 
variables, and $D$ is a Lie subalgebra of $W_n = Der(O_n)$, the simple Lie algebra of the 
general Cartan type. To ensure semisimplicity, it is assumed that $O_n$ does not contain proper $D$-invariant ideals.

\begin{theorem}
Let $n \ge 1$. Then $Z^2_{comm} ((S \otimes O_n) \inplus D) \simeq Z^2_{comm} (D)$.
\end{theorem}

\begin{proof}
According to Lemma \ref{lemma-yaya}, $Z^2_{comm} ((S \otimes O_n) \inplus D)$
is isomorphic to
\begin{align*}
& Z^2_{comm}(S) \otimes \set{\chi\in O_n^*}{\chi(d(f)g - fd(g)) = 0 \text{ for any } f,g\in O_n, d\in D} 
\\ \oplus &
C(S) \otimes \set{\beta \in HC^1(O_n)}{\beta(d(f),g) - \beta(f,d(g)) = 0 
\text{ for any } f,g\in O_n, d\in D} 
\\ \oplus & Z^2_{comm}(D) .
\end{align*}

Here again, the second tensor factor in the first direct summand vanishes.
Indeed, substitution of $g=1$ in the defining condition for $\chi\in O_n^*$, gives $\chi(d(f)) = 0$ for any 
$f\in O_n$ and $d\in D$. But then 
$$
\chi(d(f)g) = \frac{1}{2} \chi(d(f)g - fd(g)) + \frac{1}{2} \chi(d(f)g + f(d(g))) = \chi(d(fg)) = 0
$$
for any $f,g\in O_n$, $d\in D$. The space $D(O_n)O_n$ is evidently a 
$D$-invariant ideal of $O_n$, hence it coincides with the whole $O_n$, and $\chi = 0$.

And again, $C(S) = 0$ by Lemma \ref{no-cl}, so the second direct summand vanishes 
too, and the desired isomorphism follows.
\end{proof}

This is similar in spirit to the computation of the second cohomology of some
modular semisimple Lie algebras in \cite[\S 3 and Erratum and addendum]{tams}.

\section*{Acknowledgements}

Thanks are due to the anonymous referee for pointing few inaccuracies in and improvements
to the previous version of the manuscript,
to Dimitry Leites for helpful comments, and to J. Marshall Osborn who kindly 
supplied a reprint of his (hard-to-find) paper \cite{osborn}.

\end{document}